\documentclass{slides}
\usepackage{amssymb,latexsym,amsmath,amsthm,graphics,lattice}
\usepackage{graphics}
\usepackage{amscd}

\newtheorem{theorem}{Theorem}
\newtheorem{lemma}{Lemma}

\newcommand{\fg}[1]{{\mathbf M}_3\langle#1\rangle}

\newcommand{\cpe}{con\-gruence-pre\-serv\-ing extension\xspace}

\newcommand{\Ji}[1]{\tup{J}(#1)}

\newcommand{\R}{\D R}

\newcommand{\disp}[1]{\centerline{\tbf{\large #1}}}
\newcommand{\dispsub}[1]{\centerline{\tbf{#1}}}

\newtheorem*{problemn}{Problem}
\DeclareMathOperator{\Con}{Con}

\DeclareMathOperator{\Id}{Id}

\DeclareMathOperator{\rs}{rs}

\newcommand{\prodm}[2]{\prod(\,#1\mid#2\,)}

\begin{document}

\title{Representing\\ finite distributive lattices\\as congruence lattices of lattices\\
Tulane University} 
\author{G. Gr\"atzer}

\date{}
\maketitle 
 
\begin{slide}

\centerline{\includegraphics{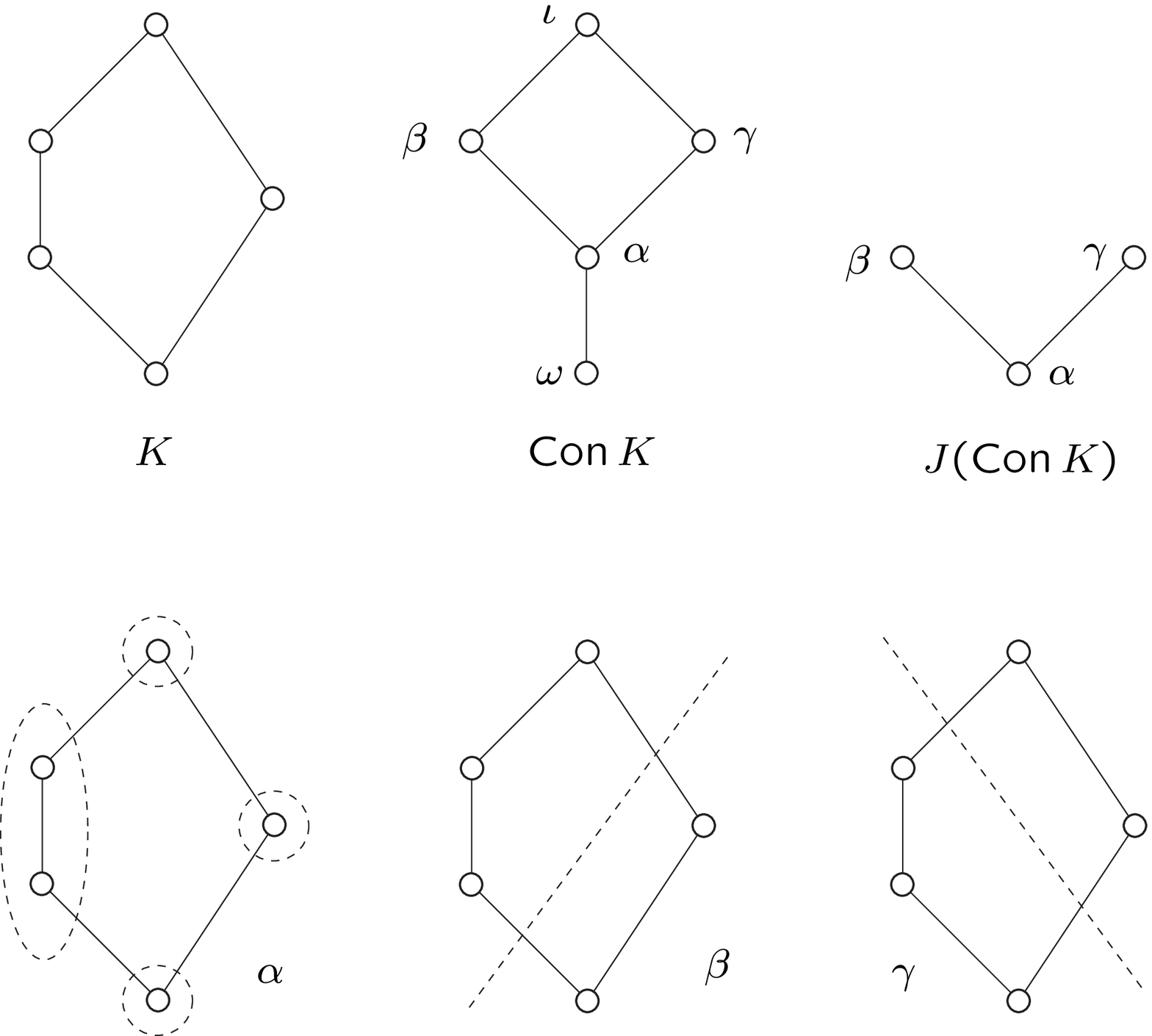}}
\end{slide} 

\begin{slide}

\centerline{\includegraphics{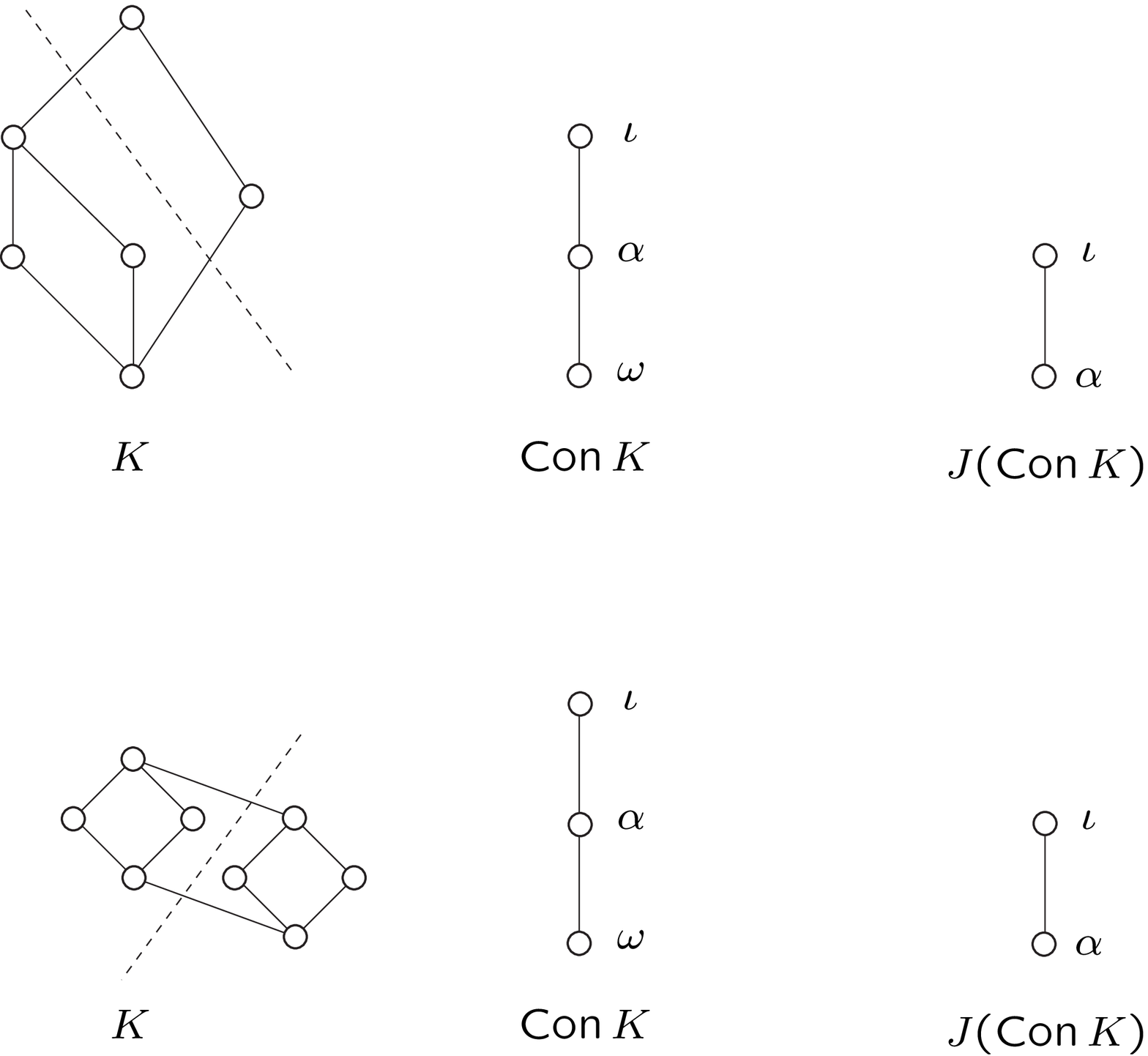}}
\end{slide} 

\begin{slide}

\vspace{50pt}

\tbf{Dilworth's theorem.}
Every finite distributive lattice $D$ can be represented as the congruence lattice of a
finite lattice $L$.

\medskip
We want:

Every finite distributive lattice $D$ can be represented as the congruence lattice of a
\tbf{nice} finite lattice $L$.

G. Gr\"atzer and E. T. Schmidt, 1962:

nice = sectionally complemented 

\end{slide}

\begin{slide}

\centerline{\includegraphics{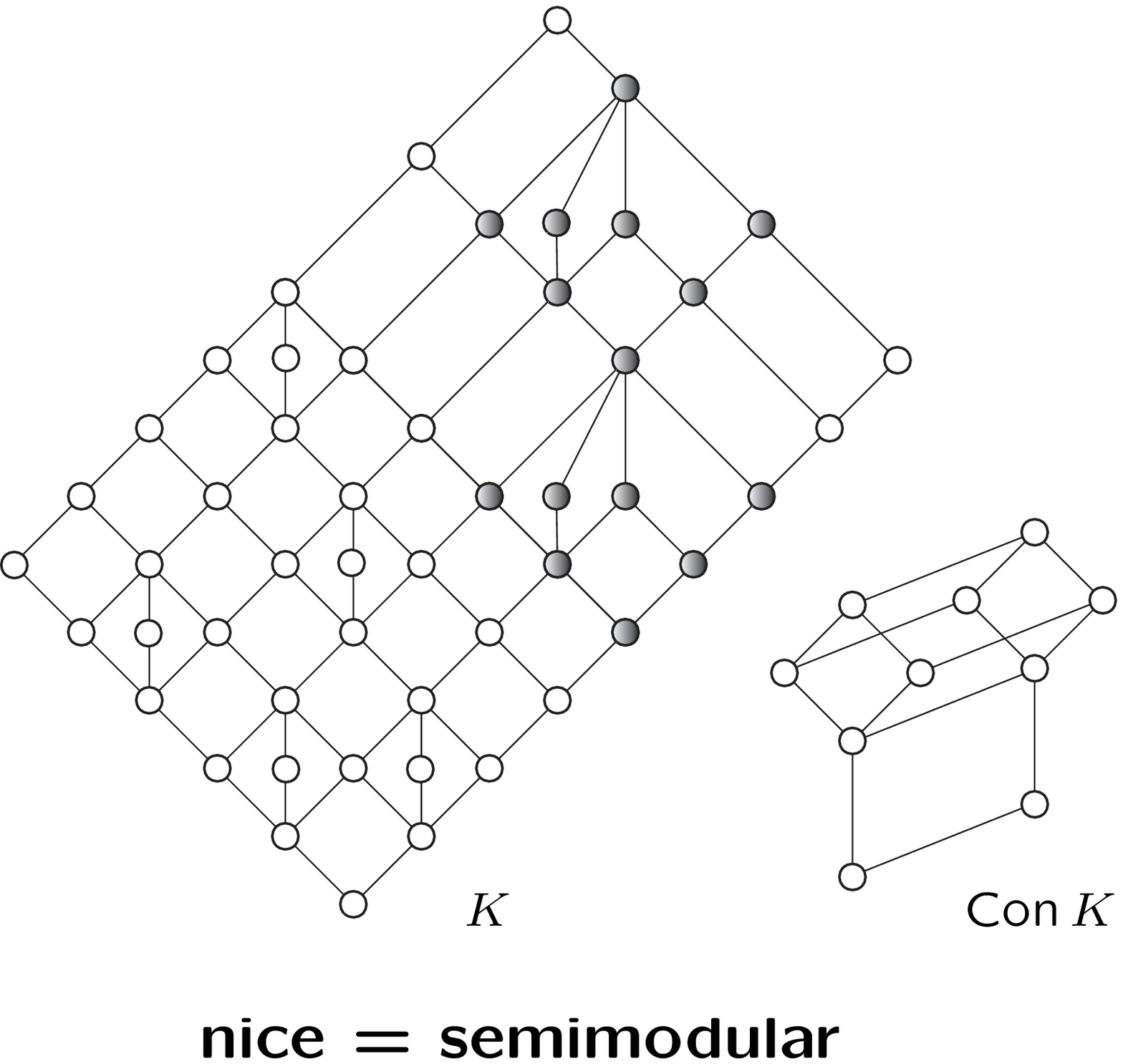}}

\end{slide} 

\begin{slide}

\centerline{\includegraphics{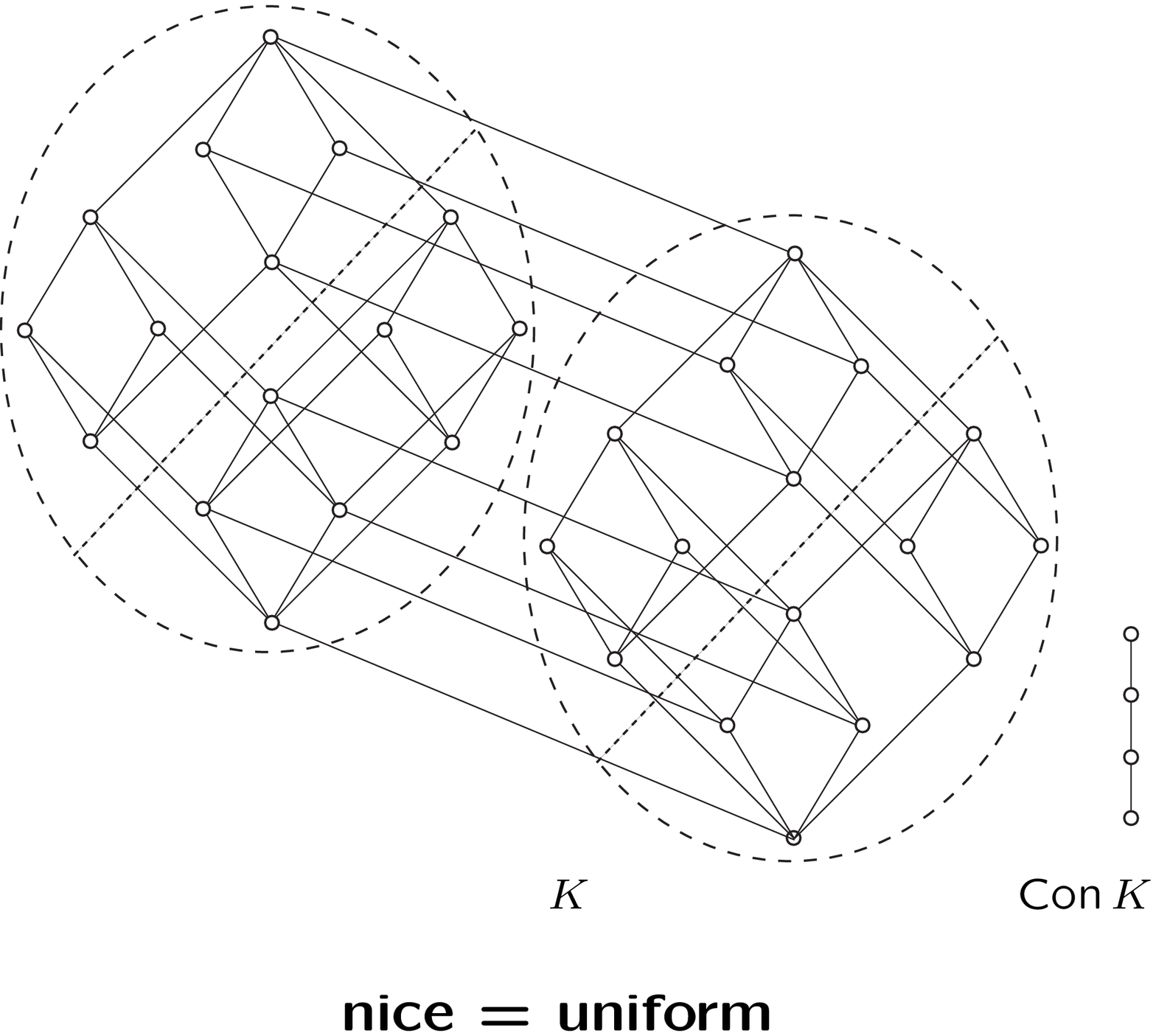}} 

\end{slide} 

\begin{slide}

\centerline{\includegraphics{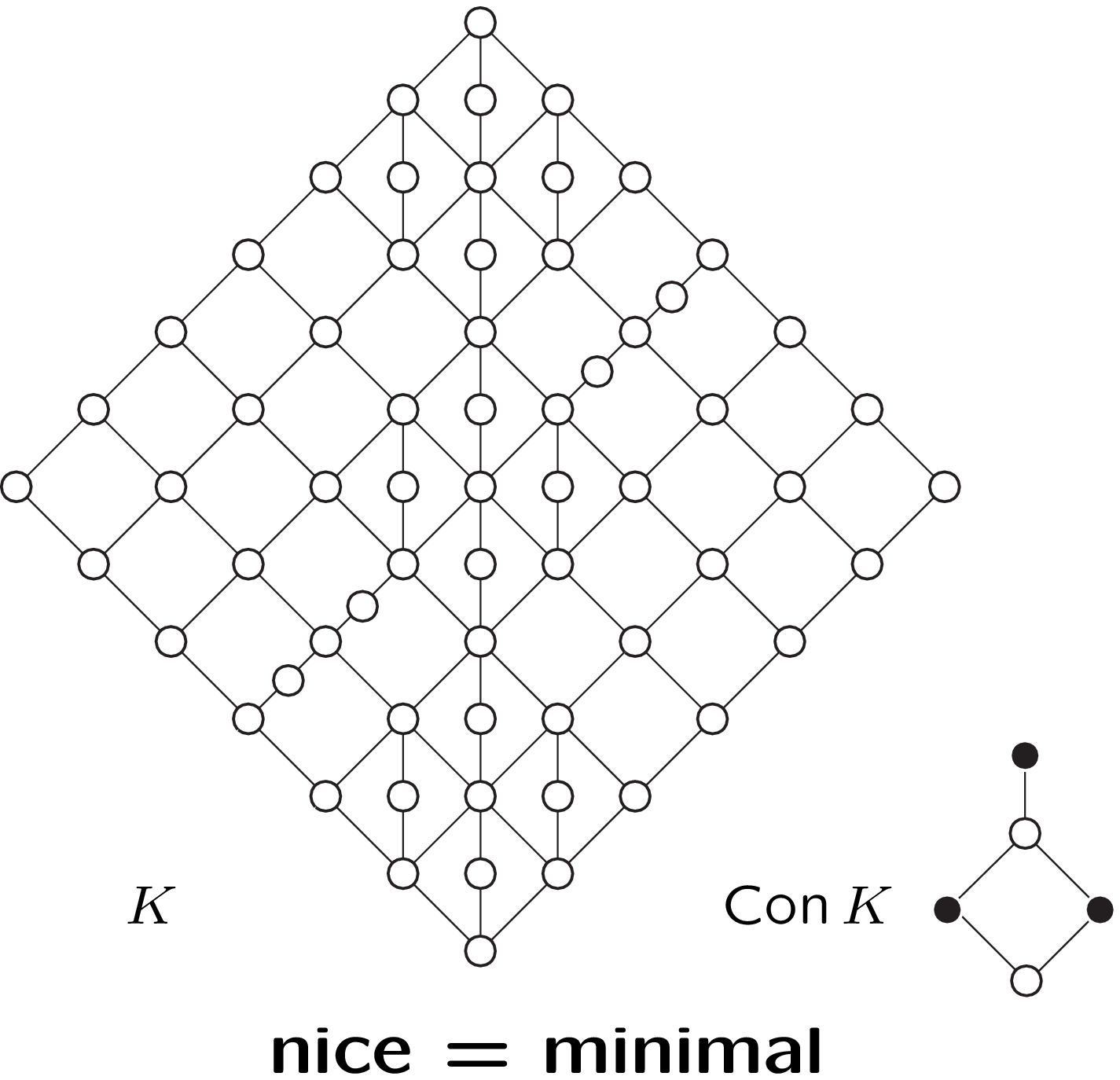}}

\end{slide} 

\begin{slide}

A \emph{finite distributive lattice} $D$ is determined by the poset $\Ji D$ of
\emph{join-irreducible elements}. So a representation of a finite distributive lattice
$D$ as the congruence lattice of a lattice $L$ is really a representation of a finite
poset $P$ ($= \Ji D$) as the poset of join-irreducible congruences of a finite lattice
$L$. 

We want:

Every finite poset $P$ can be represented as the poset of join-irreducible congruences
of a \tbf{nice} finite lattice $L$.

\end{slide}  
\begin{slide} 

There are two types of representation theorems: \\
(a) The straight representation theorems.\\
(b) The \cpe results.\\

\end{slide} 
\begin{slide}

Let $K$ be a finite lattice.  

A finite lattice $L$ is a \emph{\cpe} of $K$, if $L$ is an
extension and every congruence $\gQ$  of $K$ has \emph{exactly one} extension
$\ol \gQ$ to~$L$---that is, $\ol \gQ|_K = \gQ$.  

Of course, then the congruence lattice of $K$ is isomorphic to the congruence lattice of
$L$. 

\bigskip

\centerline{\includegraphics{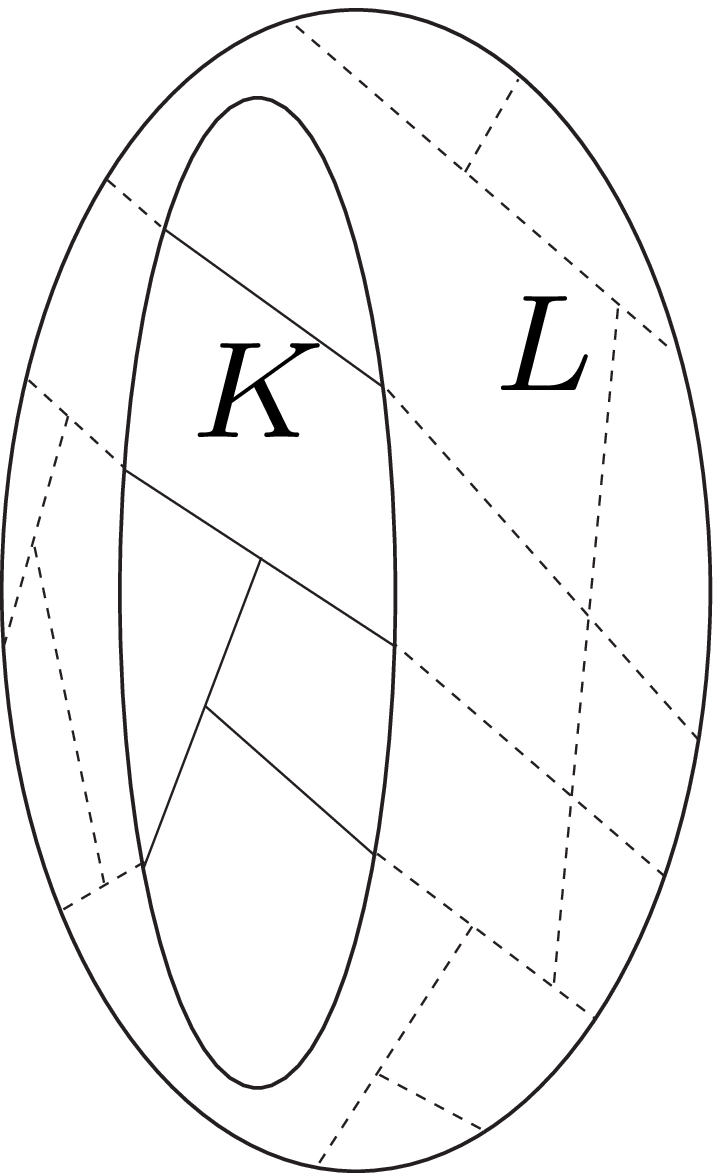}}

\end{slide} 
\begin{slide}

We could say that the congruence lattice of $L$ is ``naturally isomorphic'' to
the congruence lattice of $K$ or that the \emph{algebraic reasons} determining the
congruence lattice of $K$ are carried over to $L$. 

Examples:

\centerline{\includegraphics{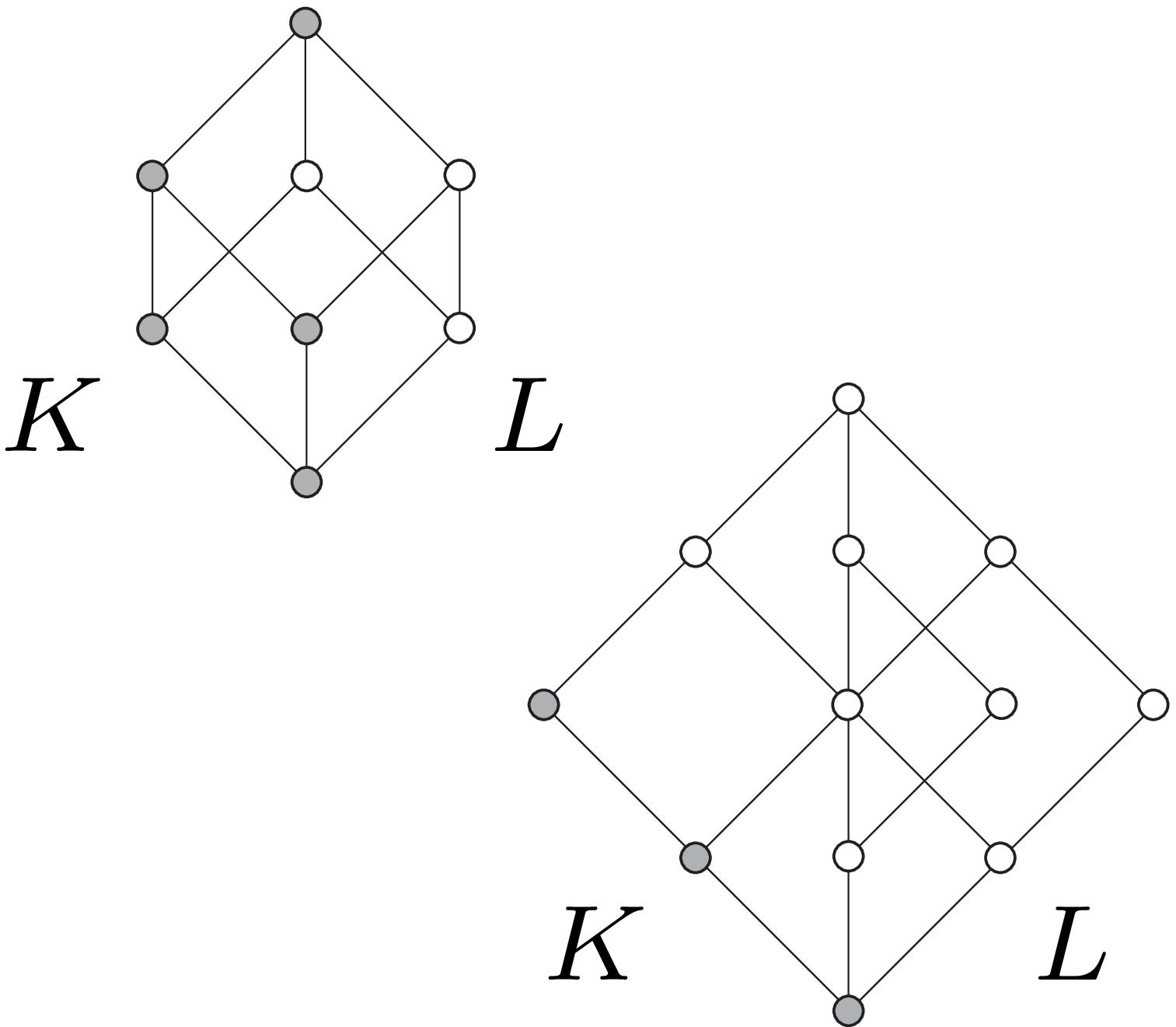}}

\end{slide} 
\begin{slide}
Not examples:

\vspace{30pt}

\centerline{\includegraphics{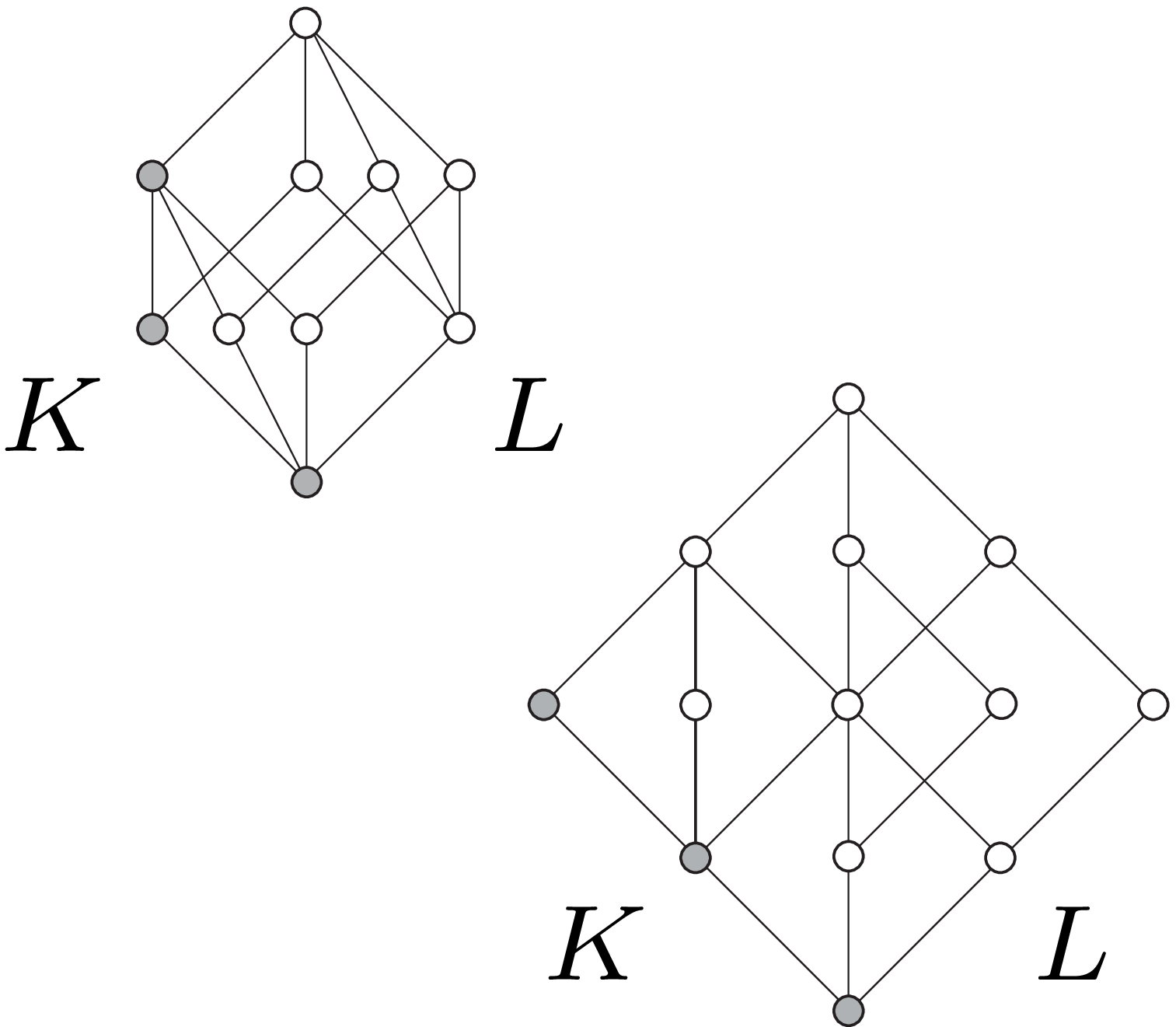}}

\end{slide} 
\begin{slide}

\disp{nice = sectionally complemented}

\vspace{30pt}

\dispsub{1. The results}

\bigskip

\begin{theorem}[G. Gr\"atzer and E.\,T. Schmidt, 1962]\label{T:relcomprepr}
Every finite distributive lattice $D$ can be represented as the congruence lattice of a
finite sectionally complemented lattice $L$.
\end{theorem}

\bigskip

\begin{theorem}[G. Gr\"atzer and E.\,T. Schmidt, 1999]\label{T:relcompcpe}
Every finite lattice $K$ has a finite, sectionally complemented, \cpe $L$.
\end{theorem}

\end{slide} 
\begin{slide}

\disp{nice = sectionally complemented}

\vspace{30pt}

\dispsub{2. Basic technique:}
\dispsub{Chopped lattices}

Let $M$ be a poset satisfying the following two conditions:\\
(i) $\inf\set{a, b}$ exists in $M$, for any $a$, $b \in M$;\\
(ii) $\sup\set{a, b}$ exists, for any $a$, $b \in M$ having a common upper bound in~$M$.

We define in $M$:
 \[
    a \mm b = \inf\set{a, b};
 \]
and
 \[
   a \jj b = \sup \set{a, b},
 \]
whenever $\sup \set{a, b}$ exists in~$M$. This makes $M$ into a
partial lattice, called a \emph{chopped lattice}. 

\end{slide} 
\begin{slide}
Examples:

\vspace{40pt}

\centerline{\includegraphics{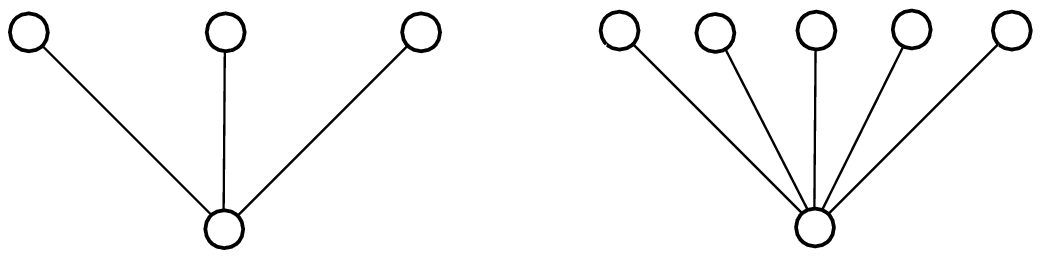}}

\vspace{40pt}

\centerline{\includegraphics{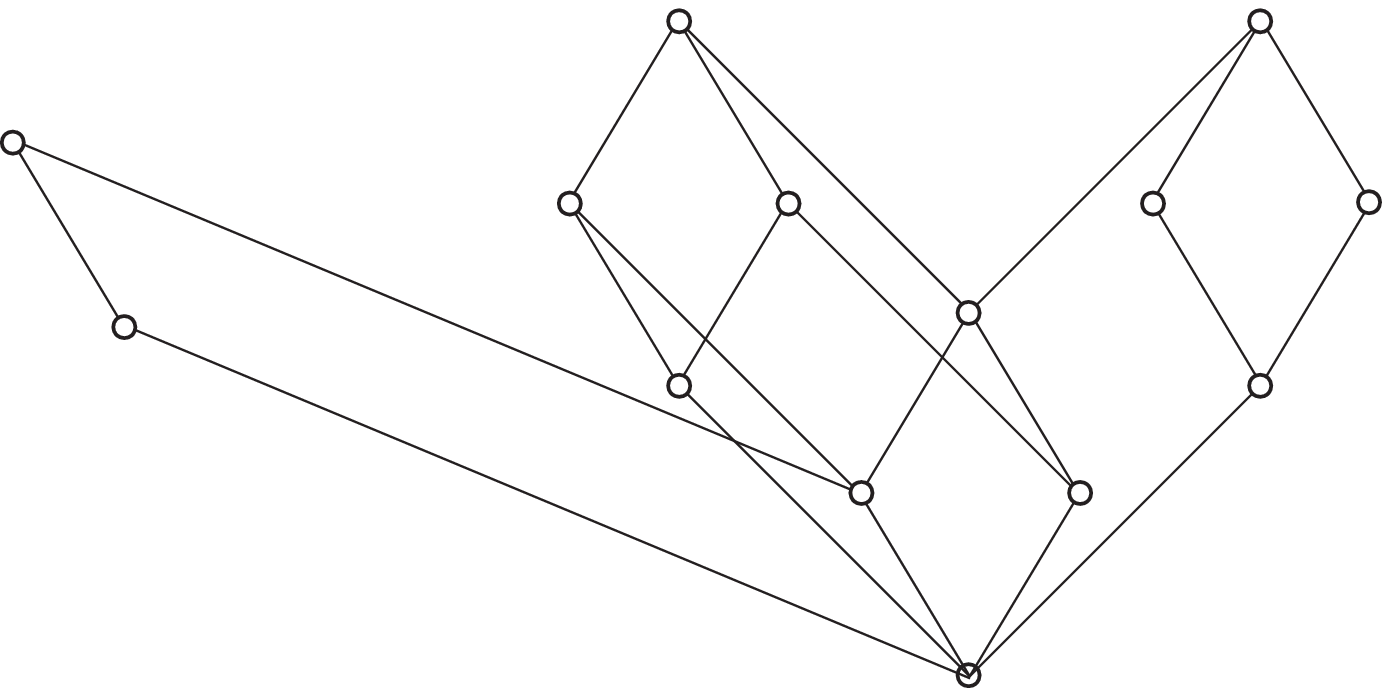}}

The second chopped lattice has a 32 element ideal lattice!

\end{slide} 
\begin{slide}

\begin{lemma}[G. Gr\"atzer and H. Lakser]\label{L:chopped} Let $M$ be a finite chopped
lattice. Then, for every congruence relation~$\gQ$, there exists exactly one congruence
relation $\ol{\gQ}$ of $\Id M$ such that, for $a$, $b \in M$,
\[
   \con {(a]} = {(b]} (\ol{\gQ}) \text{\qq if{}f\qq} \con a = b(\gQ).
\]
\end{lemma}

\end{slide} 
\begin{slide}

\begin{problemn}\label{P:seccomp}
Let $M$ be a sectionally complemented chopped lattice. Under what conditions is $\Id M$
sectionally complemented? 
\end{problemn}

\bigskip

\end{slide} 
\begin{slide}

\disp{nice = sectionally complemented}

\vspace{30pt}

\dispsub{3. The representation theorem}

We use chopped lattices to prove Theorem~\ref{T:relcomprepr}. Let $D$ be a finite
distributive lattice, and form the finite poset $P = \Ji D$.

\bigskip

\bigskip

\begin{lemma}\label{L:choppedconstruction}
Let $D$ be a finite distributive lattice. Then there exists a finite chopped lattice
$M$ such that $\Con M$ is isomorphic to~$D$.
\end{lemma}

Let us illustrate this with $D$, the four-element chain, 
\[
   0 \prec c \prec b \prec a;
\]
so $\Ji D$ is the three-element chain 
\[
   c \prec b \prec a.
\]

\end{slide} 
\begin{slide}

 Take the finite set $M_0 = \set{a, b, c} \uu \set{0}$ and make it a meet-semilattice by
defining $\inf\set{x, y} = 0$, if $x \ne y$. Note that the congruence relations of $M_0$
are in one-to-one correspondence with subsets of $\set{a, b, c}$.

\bigskip

\bigskip

\bigskip

 \centerline{\includegraphics{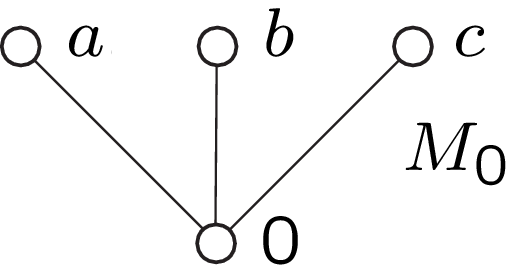}}

We must force that $a \equiv 0$ implies that $b \equiv 0$ implies that $c \equiv 0$. 

\end{slide} 
\begin{slide}

To accomplish this, we use the lattice 
\[
   N_6 = N(p, q).
\] 

\bigskip

\bigskip

\centerline{\includegraphics{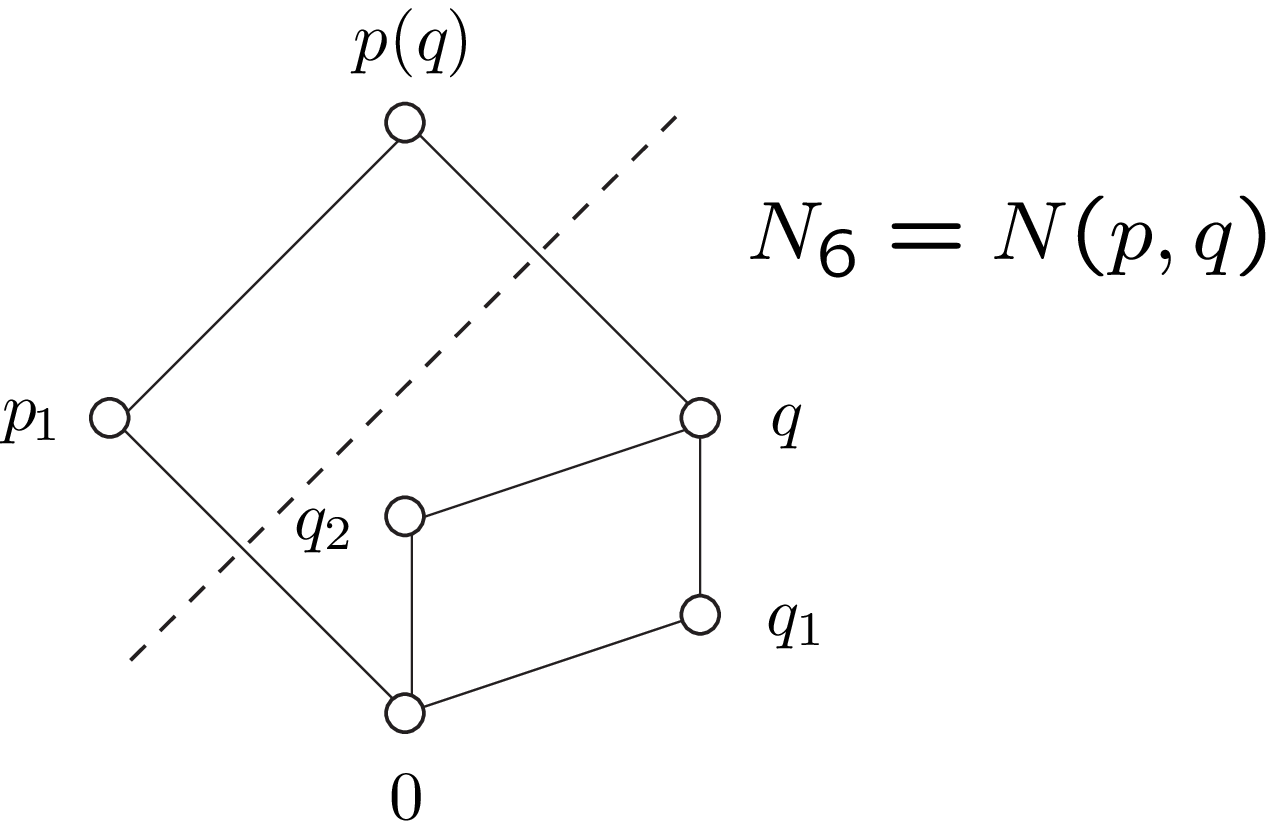}}

In $N(p, q)$, $p_1 \equiv 0$ ``implies''
that $q_1 \equiv 0$, but $q_1 \equiv 0$ ``does not imply'' that $p_1 \equiv 0$.

We construct the finite chopped lattice $M$ by ``inserting'' $N(p, q)$ in $M_0$,
for $a$, $b$ and for $b$, $c$, by appropriately doubling $b$ and $c$.

\end{slide} 
\begin{slide}

\centerline{\includegraphics{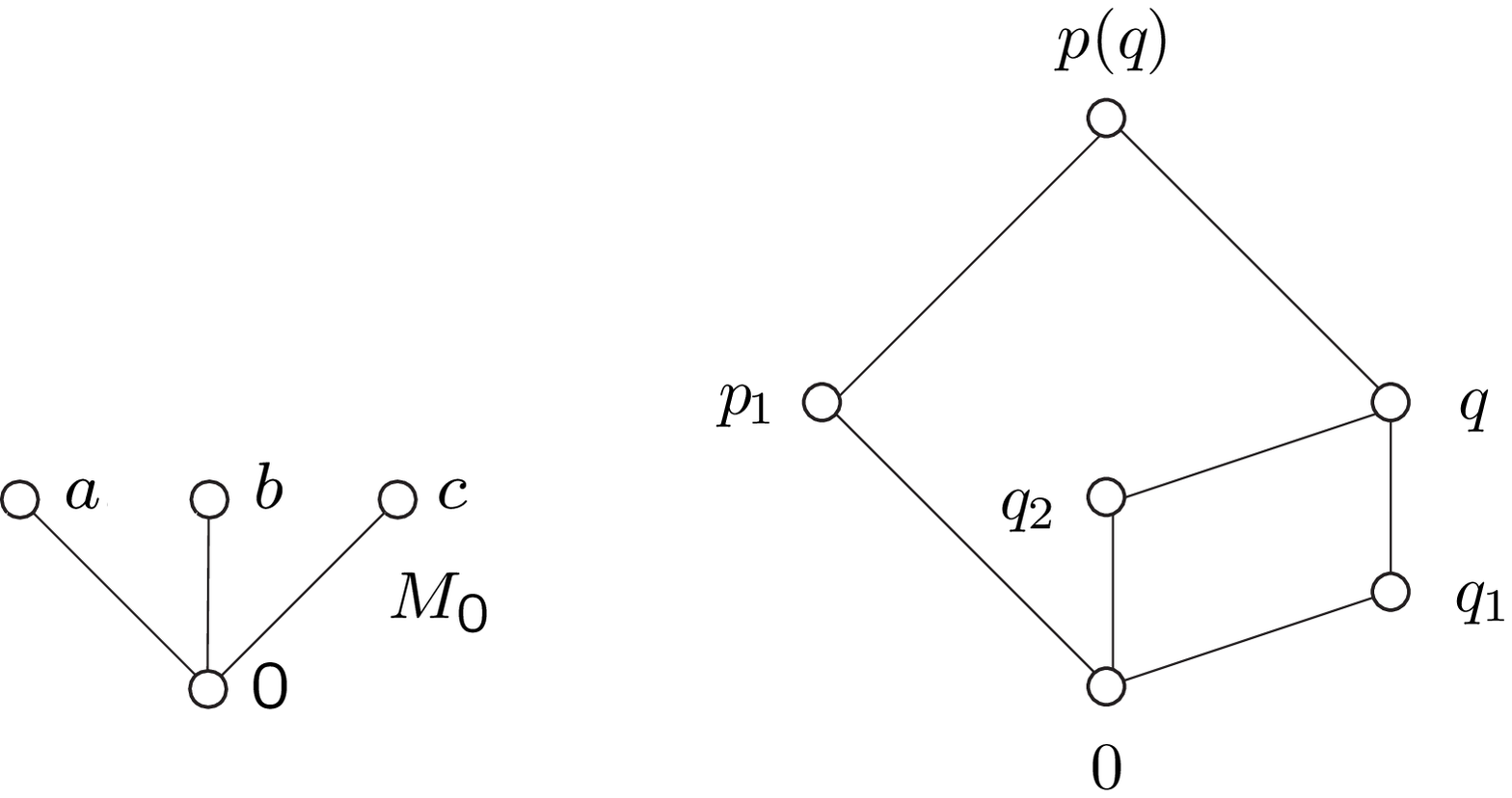}}

\bigskip

\bigskip

\centerline{\includegraphics{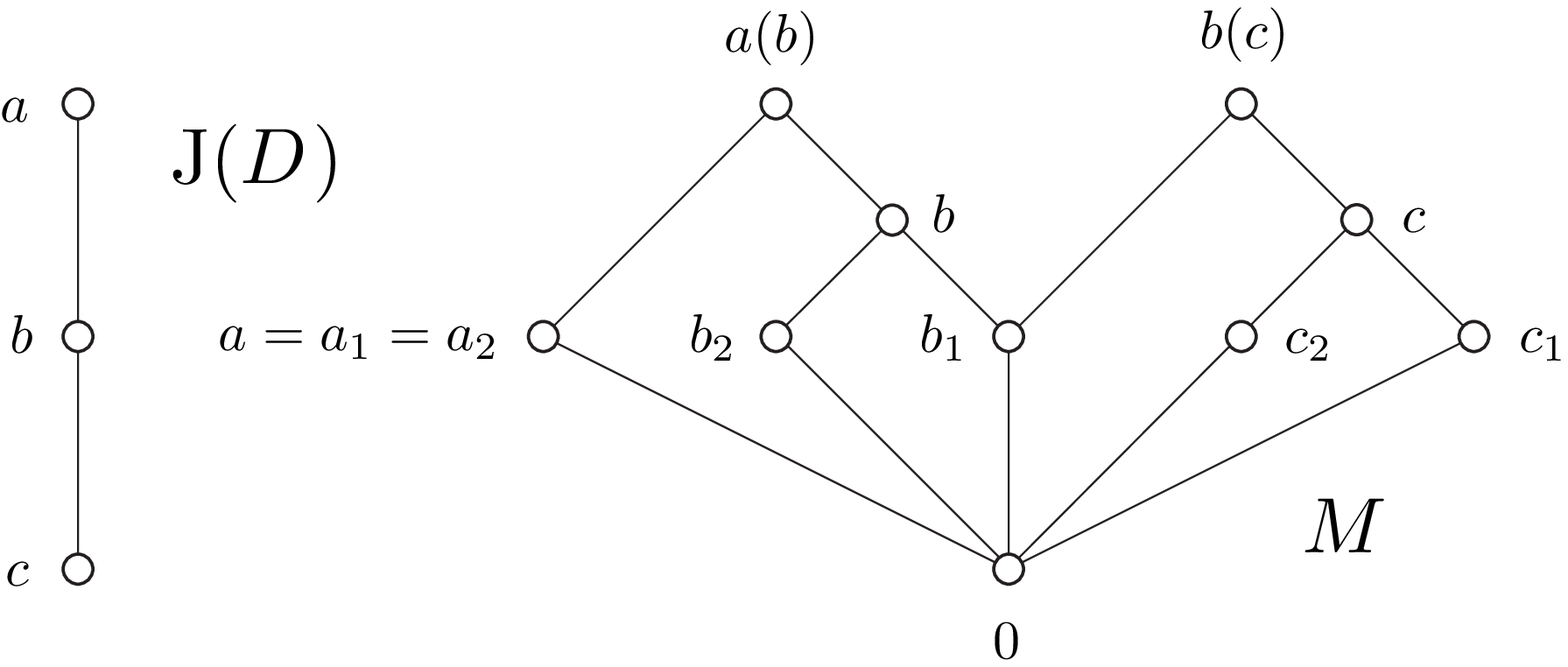}}

\end{slide}
\begin{slide}

\disp{nice = sectionally complemented}

\vspace{30pt}

\bigskip

\dispsub{4. Basic technique:}
\dispsub{Cubic extensions}

For a finite lattice $K$, let $\setm{K_i}{i \in I}$ be the subdirect factors of $K$.
For each $K_i$, we select $S(K_i)$, a finite \emph{simple} extension of $K_i$.

Now we form:
 \[
   \R(K) = \prodm{S(K_i)}{I \in I},
 \]
and call it a \emph{cubic extension} of $K$.

Then there is a one-to-one correspondence between subsets of $I$ and congruences $\gQ$ of
$\R(K)$; hence, the congruence lattice of $\R(K)$ is a finite Boolean
lattice. 

Every congruence of $K$ extends to $\R(K)$. Each $S(K_i)$ can be chosen to be
sectionally complemented. Then the cubic extension is sectionally complemented.

\end{slide} 
\begin{slide}

\disp{nice = minimal}

The lattice $L$ constructed by R.\,P. Dilworth to represent $D$ is very large, it has
$O(2^{2n})$ elements.  

\bigskip

\begin{theorem}[G. Gr\"atzer, H.
Lakser, and E. T. Schmidt]\label{T:fincong}
Let $D$ be a finite distributive lattice with $n$ join-irreducible elements.  Then there
exists a planar lattice $L$ of $O(n^2)$ elements with $\Con L \iso D$.
\end{theorem}

We illustrate the proof of this result. 

\end{slide} 
\begin{slide}

This diagram shows a distributive lattice $D$, $P = \Ji D$, and the chain $C$ we form
from $P$. The chain is of length $2|P| = 6$, and the prime intervals are ``marked'' with
elements of $P$ as shown. We call this ``marking'' \emph{coloring}.

\bigskip

\centerline{\includegraphics{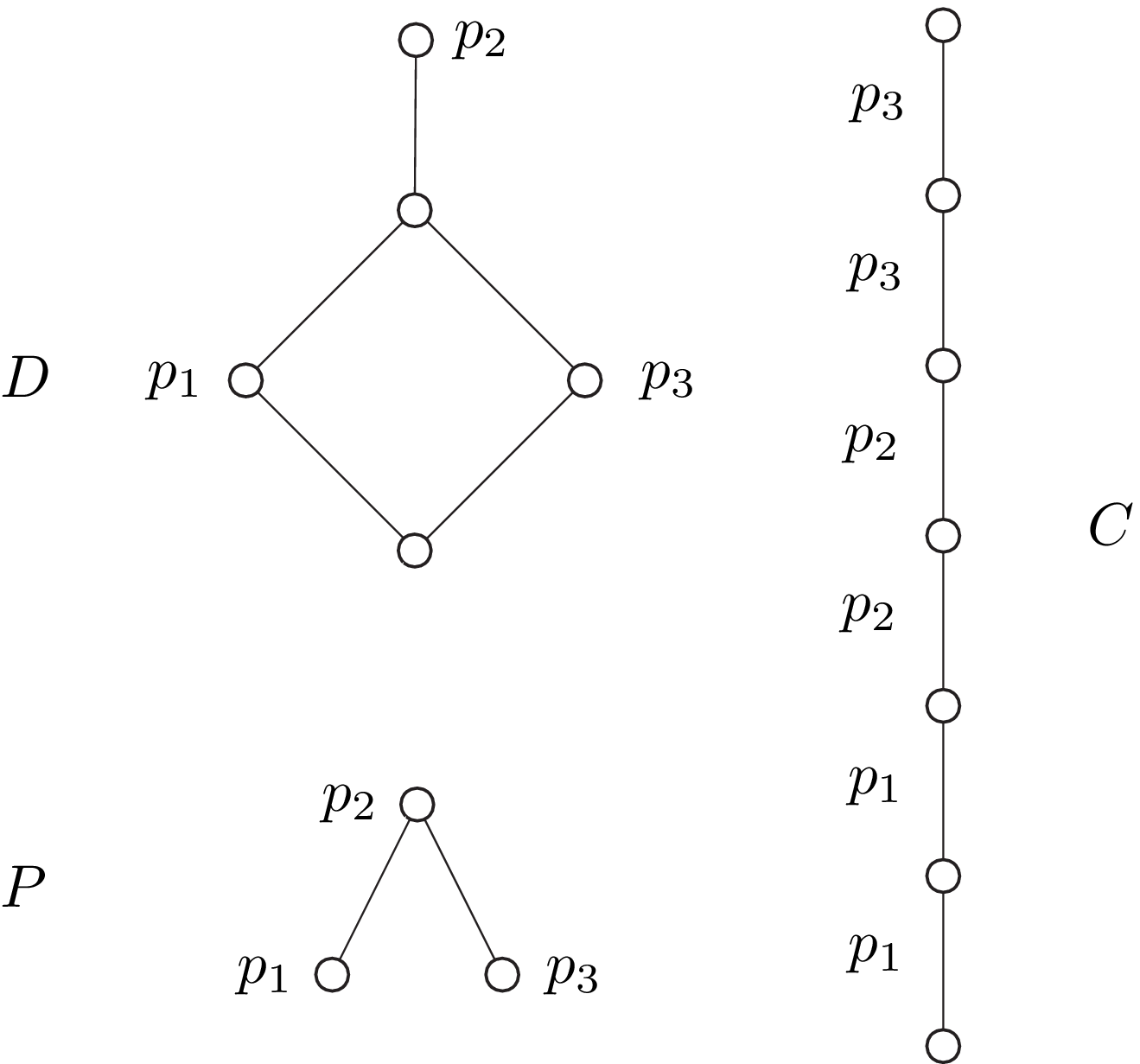}}

\end{slide} 
\begin{slide}

To construct $L$, we take $C^2$. If both lower edges of a covering square in $C^2$ have
the same color, we add an element to make it a covering $M_3$. If in $C^2$ we have a
covering $C_2 \times C_3$, where the lower $C_2$ is colored by $p$, the lower $C_3$ is
colored by $q$ twice, where $p > q$, then we add an element to make it an $N_{5, 5}$.

\bigskip

\centerline{\includegraphics{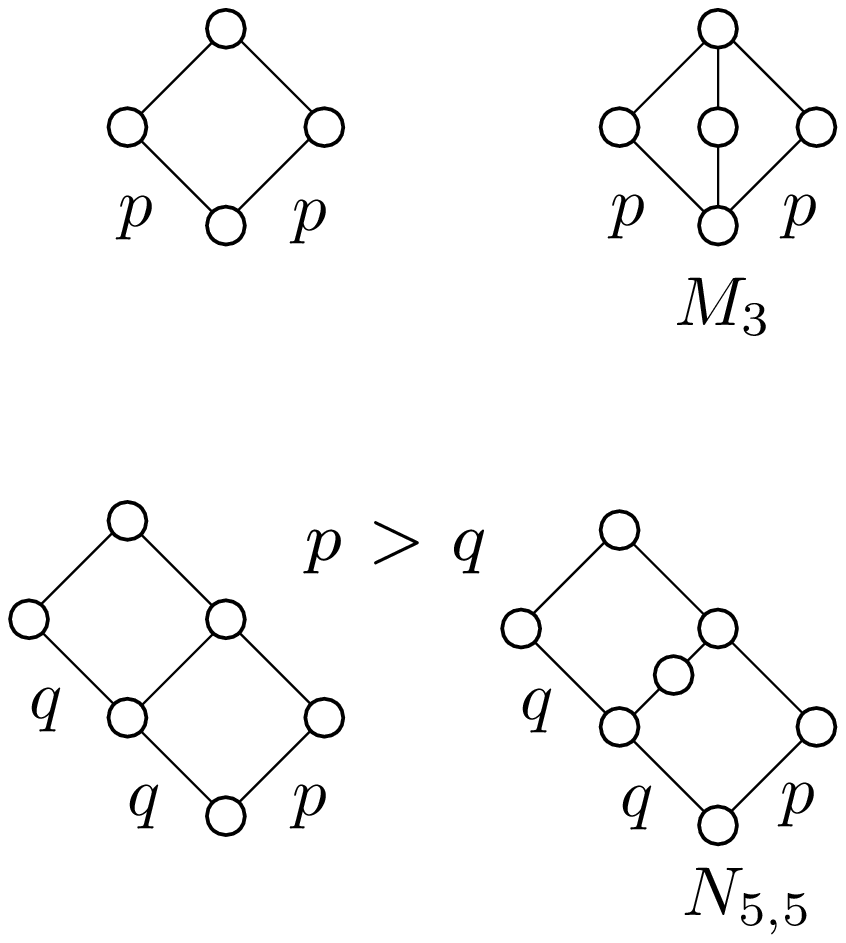}}

\end{slide} 
\begin{slide}

\centerline{Here is what we obtain:}

\bigskip

\centerline{\includegraphics{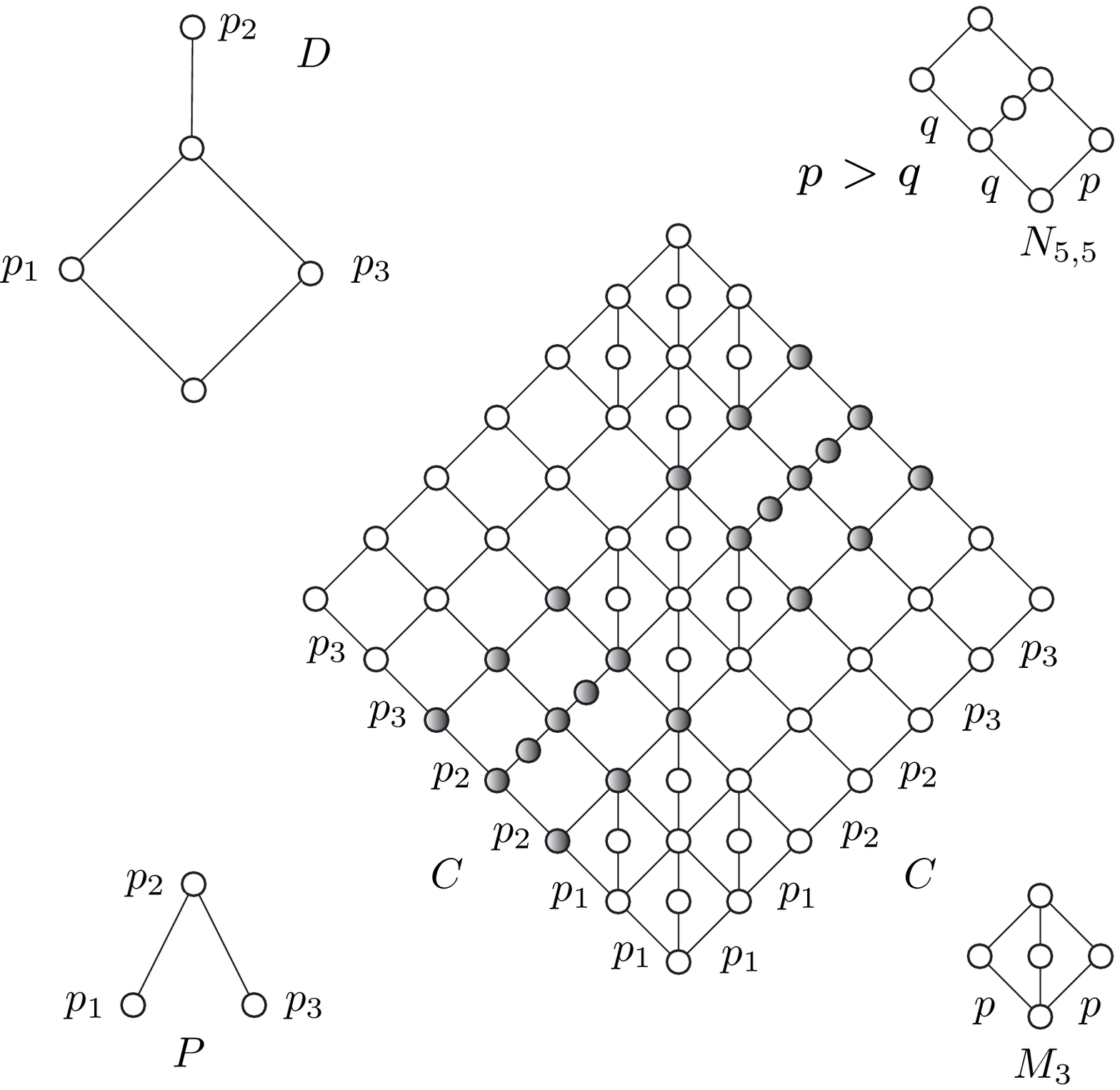}}

\end{slide} 
\begin{slide}

\begin{theorem}[G.~Gr\"atzer, I. Rival, and  N.~Zaguia, 1995]\label{T:best}
Let $\ga$ be a real number satisfying the following condition:  Every distributive
lattice $D$ with $n$ join-irreducible elements can be represented as the congruence
lattice of a lattice $L$ with $O(n^{\ga})$ elements. Then $\ga \geq 2$.
\end{theorem} 

\end{slide} 

\begin{slide}

\disp{nice = semimodular}

\vspace{30pt}

\dispsub{1. Representation theorem}

\bigskip

\bigskip

\begin{theorem}[G.~Gr\"a\-tzer, H. Lakser, and E.\,T. Schmidt, 1998]\label{T:semimod}
Every finite distributive lattice $D$ can be represented as the congruence lattice of a
finite \emph{semimodular} lattice~$S$.  In fact, $S$ can be constructed as a
\emph{planar} lattice of size $O(n^3)$, where $n$ is the number of join-irreducible
elements of $D$.
\end{theorem}

\end{slide} 
\begin{slide}

 The proof of this result is very similar to the proof of Theorem~\ref{T:fincong}. The
basic building block is $S_8$; we show two views of this lattice:

\bigskip

\bigskip

\bigskip

\centerline{\includegraphics{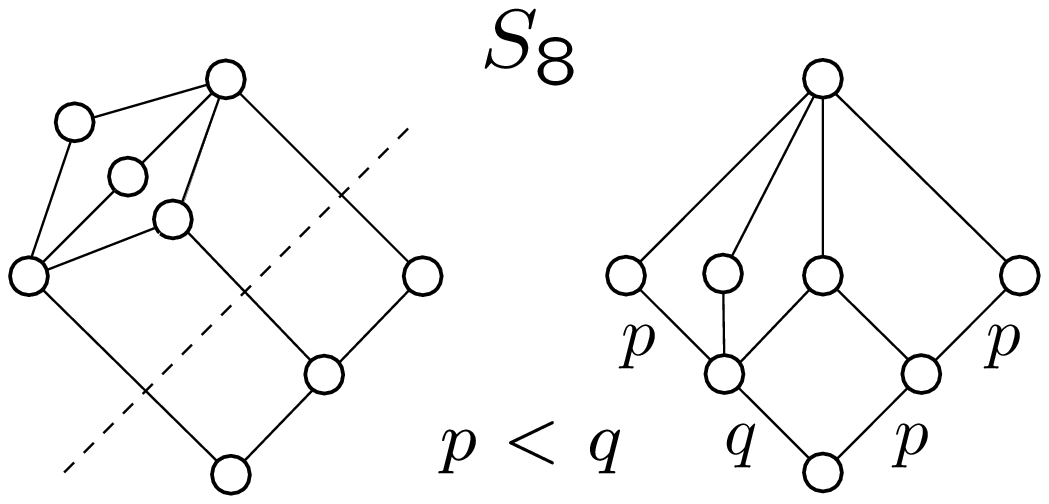}}

\bigskip

\bigskip

The first view shows $S_8$ with its only nontrivial congruence indicated with a dashed
line; the second views shows $S_8$ as it is used in the construction.

\end{slide} 
\begin{slide}

To illustrate the construction, take the following distributive lattice $D$; the poset
$\Ji{D}$ is also shown.

\bigskip

\bigskip

\bigskip

\centerline{\includegraphics{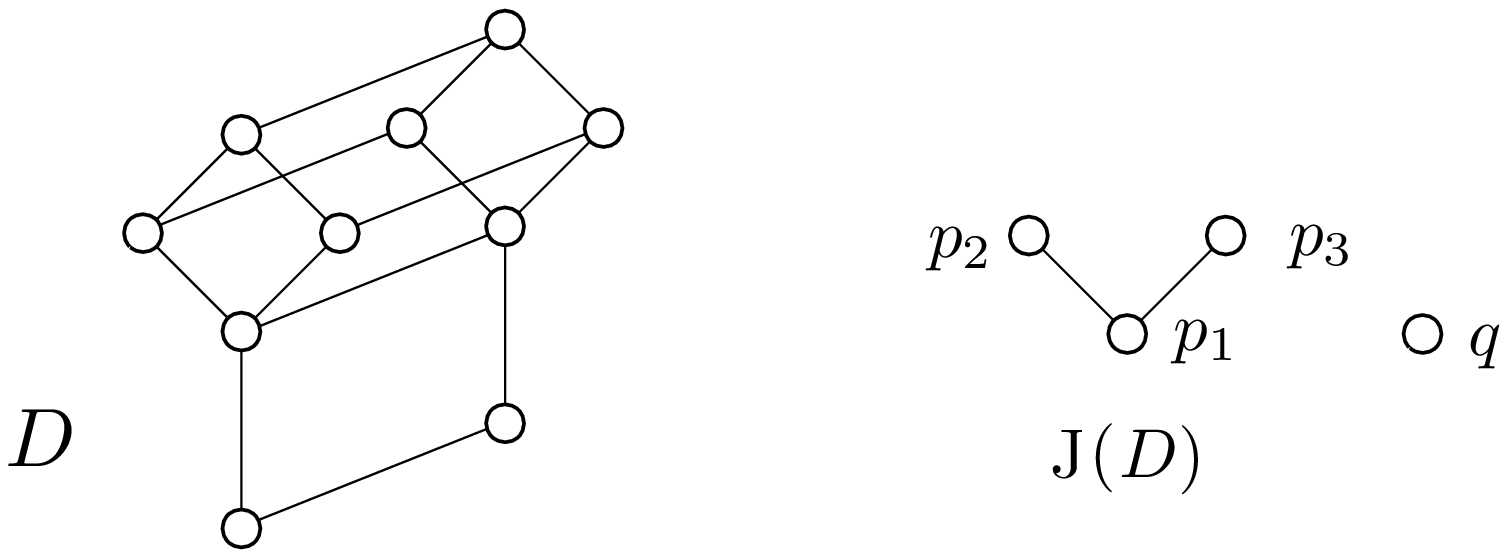}}

Now we form two chains and their direct products. We augment a covering $C_2 \times C_2$
colored by the same element on both sides by an element to form $M_3$.
For 
$p_1 \prec p_2$, we replace the covering $C_3 \times C_3$ colored by $p_2$, $p_1$ and
$p_1$, $p_1$ by $S_8$. We make a similar replacement for $p_1 \prec p_3$, to obtain $S$:

\end{slide} 
\begin{slide}

\bigskip

\centerline{\includegraphics{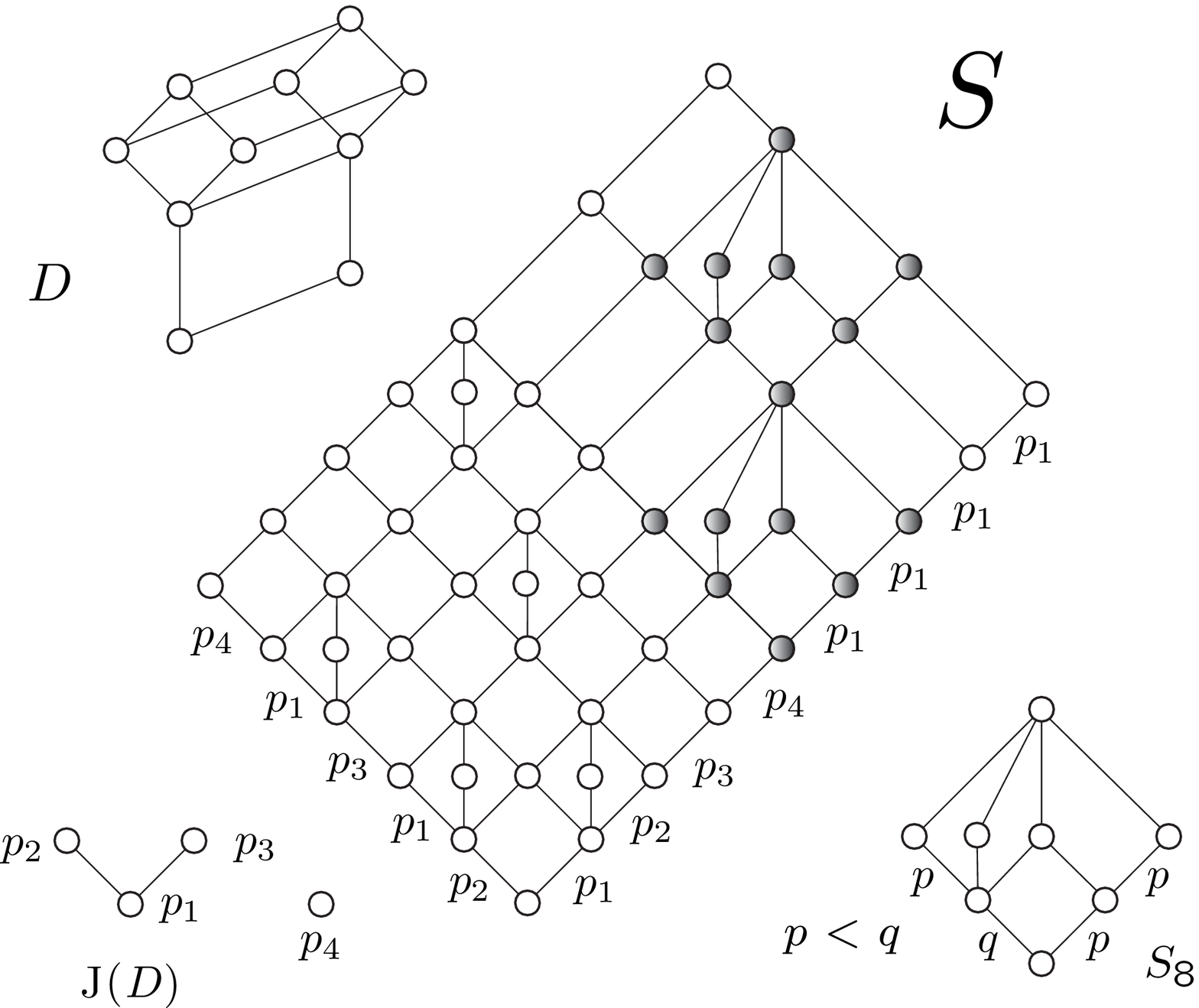}}

\end{slide} 
\begin{slide}

\disp{nice = semimodular}

\vspace{30pt}

\dispsub{2. Constructing a}
\dispsub{congruence-preserving}
\dispsub{semimodular extension}

\bigskip

\bigskip

\begin{theorem}[G. Gr\"atzer and E.\,T. Schmidt, 2001]\label{T:semimodcpe}
 Every finite lattice $K$ has a congruence-preserving embedding into a finite
\emph{semimodular} lattice $L$. 
\end{theorem}

The proof starts out with the cubic extension $\R(K)$ of $K$, where we choose each  
$S(K_i)$ semimodular. So the cubic extension is semimodular. The
congruences then are represented in a dual ideal $F$ of $\R(K)$ that is
Boolean. By~gluing a suitable modular lattice $M$ to $\R(K)$. The congruences are then
represented on a dual ideal $E'$ of $M$ that is a chain, so the proof is completed by
gluing the lattice $S$ to the construct:

\end{slide} 
\begin{slide}

\centerline{\includegraphics{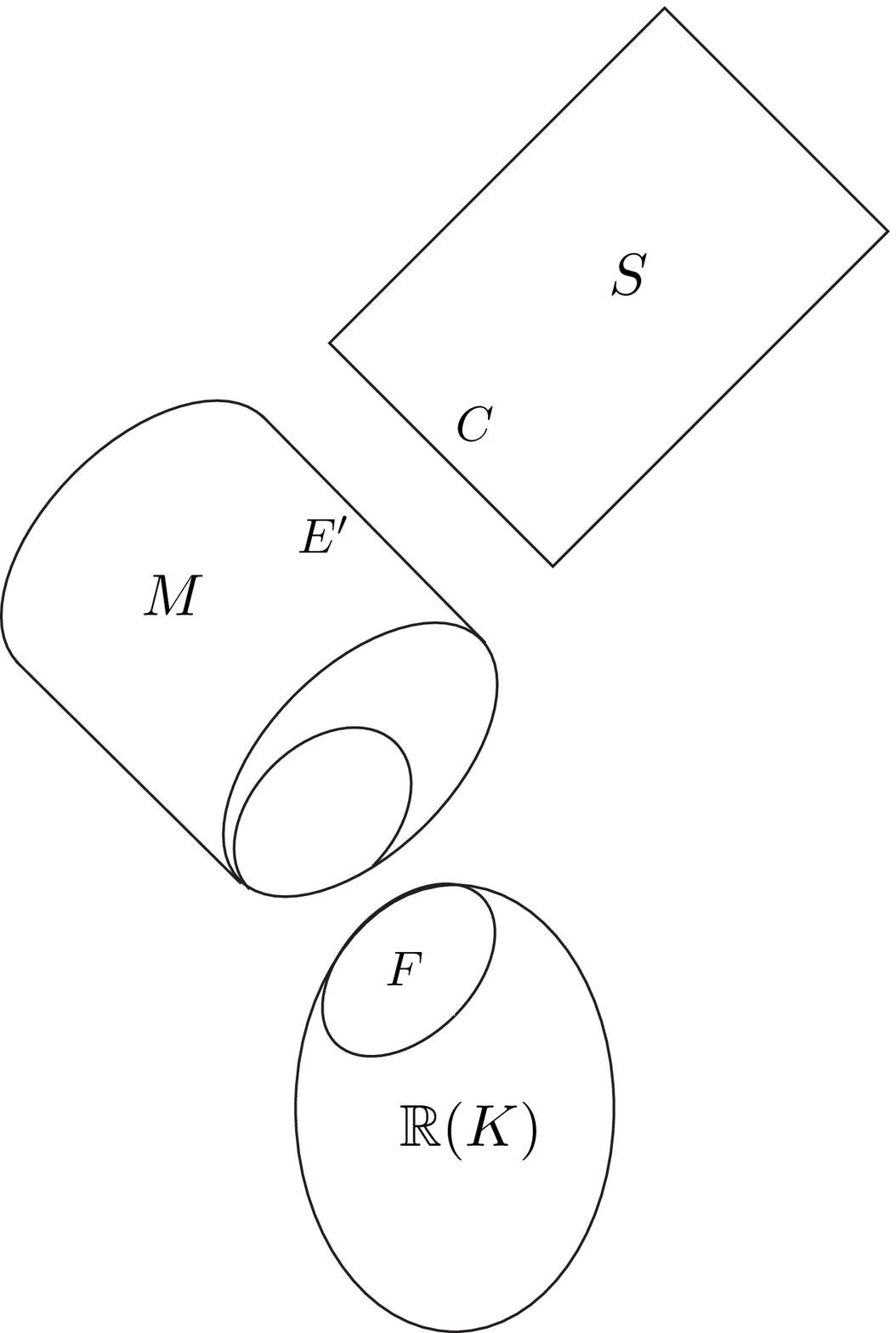}}

\end{slide} 
\begin{slide}

\disp{nice = given automorphism group}

\bigskip

\bigskip

\begin{theorem}[The Independence Theorem, V.~A. Baranski\u\i\ and
A.~Urqu\-hart, 1979]\label{T:IndThm} Let $D$ be a finite distributive lattice with more
than one element, and let $G$ be a finite group.  Then there exists a finite lattice $L$
such that the congruence lattice of $L$ is isomorphic to   $D$ and the automorphism
group of $L$ is isomorphic to $G$. 
\end{theorem}

This is a representation theorem. There is also a \cpe variant for
this result. 

\end{slide} 
\begin{slide}

\begin{theorem}[The Strong Independence Theorem, G.~Gr\"atzer and E.\,T.
Schmidt, 1995]\label{T:StrongHalf} Let $K$ be a finite lattice with more than one
element and let $G$ be a finite group. Then $K$ has a \cpe $L$ whose automorphism group
is isomorphic to $G$.
\end{theorem}

\end{slide} 
\begin{slide}

\disp{nice = regular}

\vspace{30pt}

\dispsub{1. The result}

Let $L$ be a lattice. We call a congruence relation $\gQ$ of $L$ \emph{regular}, if any
congruence class of $\gQ$ determines the congruence. Let us call the lattice $L$
\emph{regular}, if all congruences of $L$ are regular. 

Sectionally complemented lattices are regular, so we already have a representation
theorem. We also have a \cpe version:

\bigskip

\bigskip

\begin{theorem}[G.~Gr\"atzer and E.\,T. Schmidt, 2001]\label{T:regular}
Every finite lattice $L$ has a congruence-preserving embedding
into a finite regular lattice $\tilde{L}$. 
\end{theorem}

\end{slide}
 \begin{slide}

\disp{nice = regular}

\vspace{30pt}

\dispsub{2. Basic technique:}
\dispsub{Boolean triples}

For a bounded lattice $K$, let us call the triple $\vv<x,y,z> \in L^3$ \emph{boolean}
if{}f the following equations hold:
\begin{align} 
    x &= (x \jj y) \mm (x \jj z),\notag\\ 
    y &= (y \jj x) \mm (y \jj z),\notag\\ z &= (z \jj x) \mm (z \jj
y)\notag. 
\end{align}
We denote by $\fg{K} \ci K^3$ the poset of all boolean triples of
$K$. $\fg{K}$ is a bounded lattice. We identify the lattice $K$ with the
interval $[\vv<0, 0, 0>, \vv<1, 0, 0>]$ under the isomorphism $x \mapsto \vv<x, 0, 0>$.

\end{slide} 
\begin{slide}
G. Gr\"atzer and F. Wehrung, 1997:

\bigskip

\begin{theorem}\label{T:booleantriples}
$\fg{K}$ is a \cpe of $K$. 
\end{theorem}

\end{slide} 
\begin{slide}

\disp{nice = uniform}

Let $L$ be a lattice. We call a congruence relation $\gQ$ of $L$ \emph{uniform}, if any
two congruence classes of $\gQ$ are of the same size (cardinality). Let us call the
lattice $L$ \emph{uniform}, if all congruences of $L$ are uniform. 

\bigskip

\begin{theorem}[G.~Gr\"atzer, E.\,T. Schmidt, and K. Thomsen, 2002]\label{T:Uniform}
Every finite distributive lattice $D$ can be represented as the congruence lattice of a
finite uniform lattice~$L$. 
\end{theorem}

A uniform lattice is always regular, so the lattice $L$ of this theorem
is also regular.

\end{slide} 
\begin{slide}

\centerline{\includegraphics{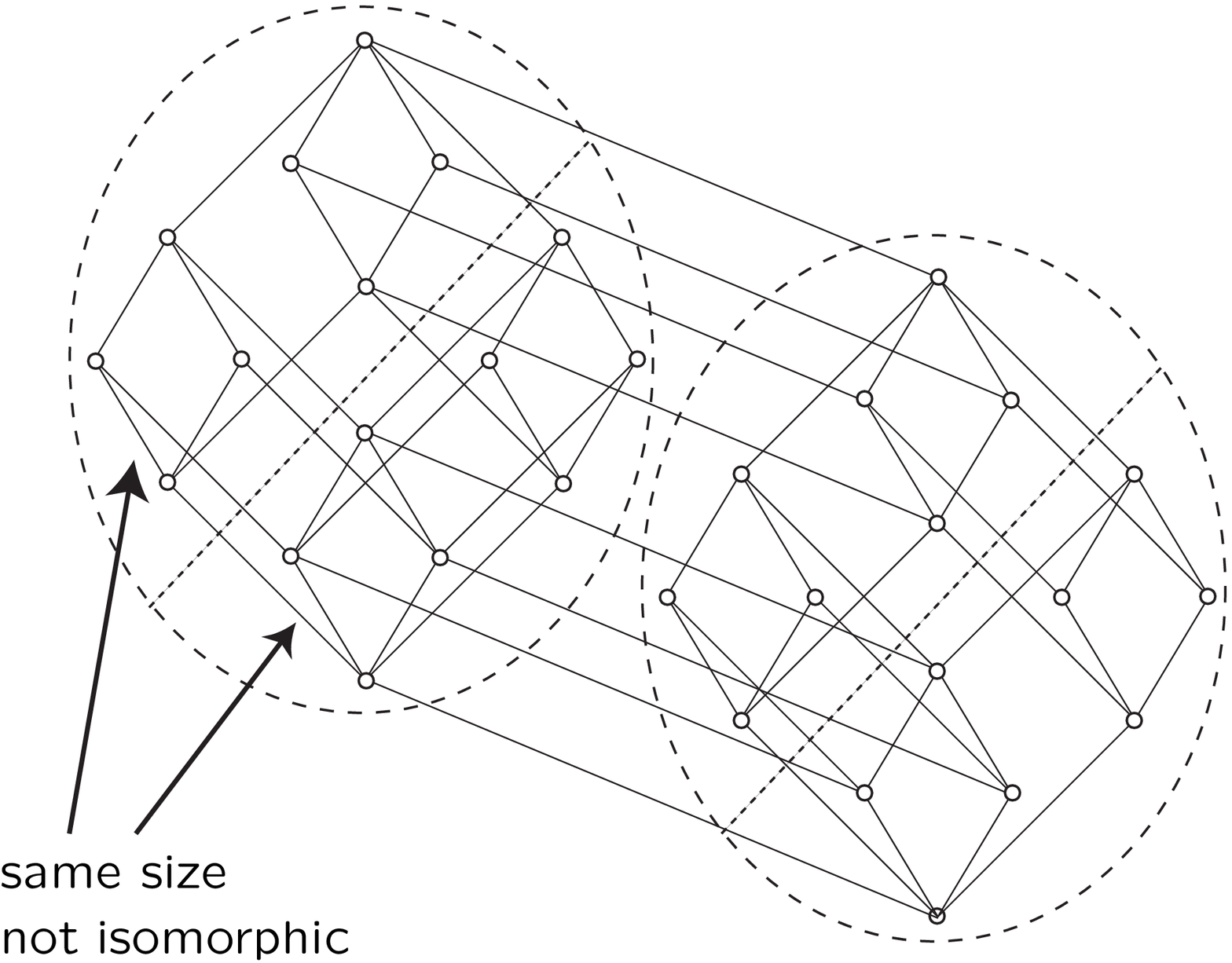}} 
The uniform construction for the four-element chain.
\end{slide} 

\begin{slide}
Why not draw less trivial examples?

If we start with the lattice:

\vspace{20pt}

\centerline{\includegraphics{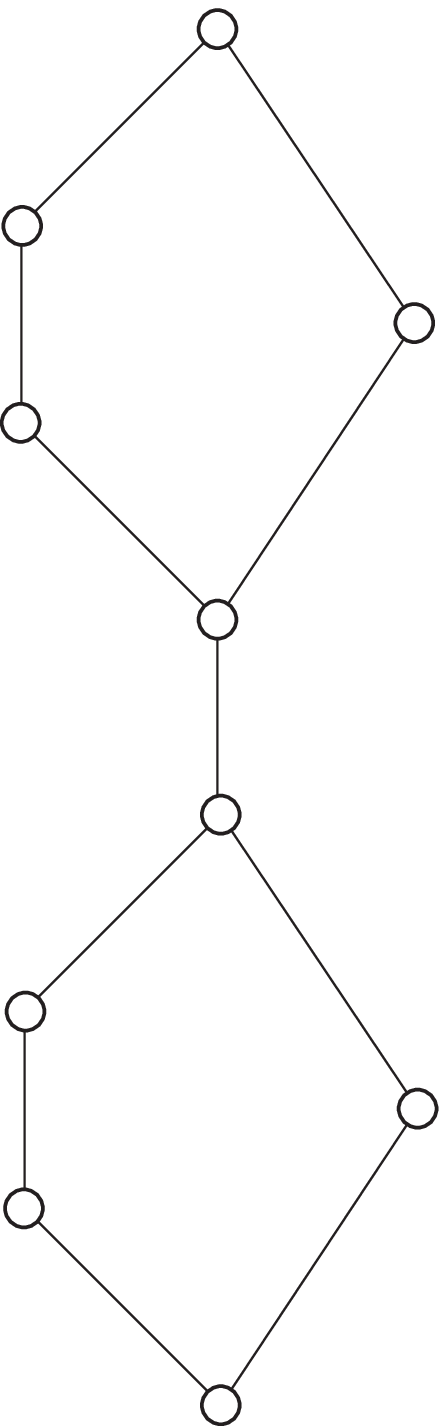}} 

then the uniform lattice we conctruct has 153,125 elements (and order dimension 7).
\end{slide} 

\begin{slide}
\disp{nice = isoform}

Let $L$ be a lattice. We call a congruence relation $\gQ$ of $L$ \emph{isoform}, if any
two congruence classes of $\gQ$ are isomorphic (as lattices). Let us call the lattice
$L$ \emph{isoform}, if all congruences of $L$ are isoform.

\smallskip

\begin{theorem}[G.~Gr\"atzer and E.\,T. Schmidt, 2002]\label{T:Isoform}
Every finite distributive lattice $D$ can be represented as the congruence lattice of a
finite, isoform lattice $L$. 

\vspace{-14pt}

\end{theorem}
Since isomorphic lattices are of the same size, this theorem is a
stronger version of the theorem for uniform lattices. 

\smallskip

\begin{theorem}[G.~Gr\"atzer, E.\,T. Schmidt, and R. W. Quackenbush, 2004]\label{T:main}
Every finite lattice $K$ has a \cpe to a finite isoform lattice $L$.
\end{theorem}

\end{slide} 

\begin{slide}

\disp{Simultaneous representations}
\disp{of two distributive lattices}

Let $L$ be a lattice and let $K$ be a sublattice of $L$. Then the \emph{restriction map}
\[
   \rs \colon \Con{L}\to\Con{K}
\]
is a $\set{0, 1, \mm}$-homomorphism. If $K$ is an ideal, then $\rs$ is a
$\set{0, 1}$-homomorphism.

G. Gr\"atzer and H. Lakser, 1986:

\bigskip

\begin{theorem}
Let $D$ and $E$ be finite distributive lattices; let $D$ have more than one element. Let
$\gf$ be a $\set{0, 1}$-homomorphism of $D$ into $E$.  Then there exists
a (sectionally complemented) finite lattice $L$ and an ideal $K$ of $L$ such that $D \iso
\Con L$,
$E
\iso
\Con K$, and $\gf$ is represented by $\rs$, the restriction map.
\end{theorem}

\end{slide} 
\begin{slide}

Example:

\centerline{\includegraphics{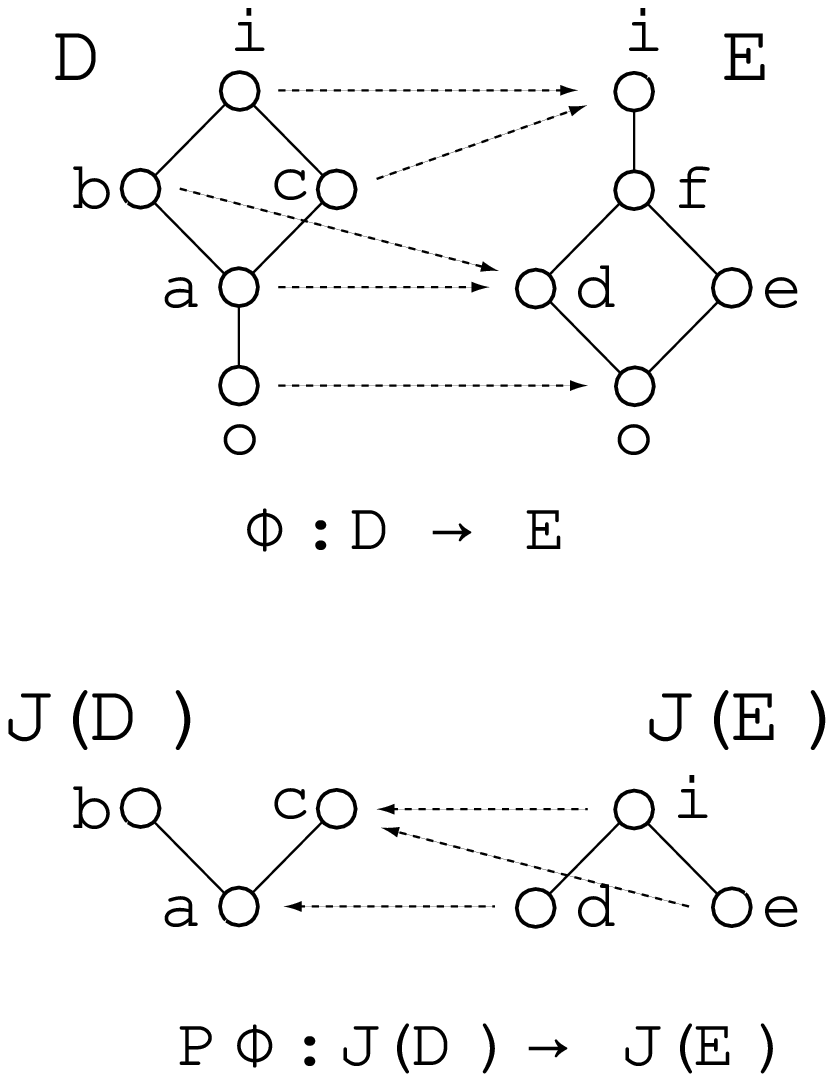}}

\end{slide} 
\begin{slide}

Building block:

\bigskip

\centerline{\includegraphics{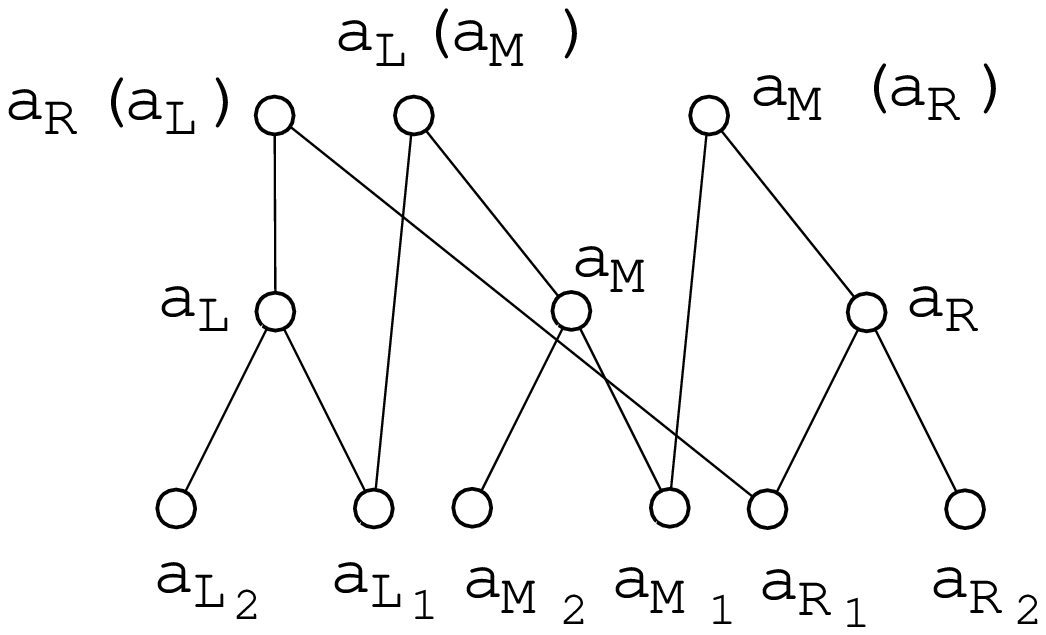}}

We take six: for $o$, $b$, $c$, $d$, $e$, $i$.

\end{slide} 
\begin{slide}

We code $\Ji D$ and $\Ji E$ with the middle elements (the $M$-s), for example, $a \prec
b$ is coded by

\bigskip

\centerline{\includegraphics{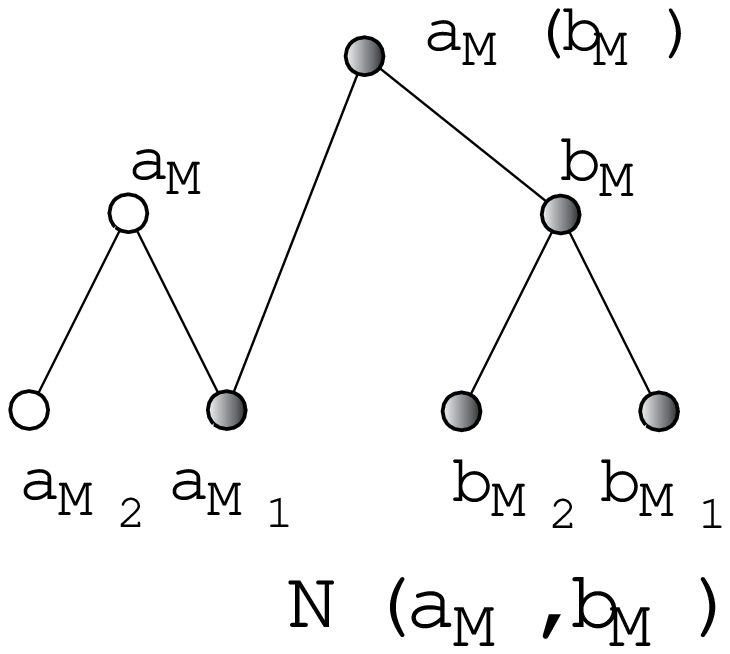}}

$P\gf$ is an equivalence: one direction we code with the $L$-s, the other, with the
$R$-s.
\end{slide} 

\begin{slide}
New problems: Make $K$ and $L$ ``nice''. 

\end{slide}

\begin{slide}
Sectionally complemented lattices: 

G. Gr\"atzer and H. Lakser, 
\emph{Notes on sectionally complemented lattices. I. 
Characterizing the 1960 sectional complement.}
Acta Math. Hungar. \tbf{108} (2005), 115--125.

G. Gr\"atzer and H. Lakser, 
\emph{Notes on sectionally complemented lattices. II. Generalizing the 
1960 sectional complement with an application to congruence
restrictions.}
Acta Math.\\ Hungar. \tbf{108} (2005), 251--258.

G. Gr\"atzer and E.\,T. Schmidt, 
\emph{Finite lattices with isoform congruences.}
Tatra Mt. Math. Publ. \tbf{27} (2003), 111--124.

Isoform lattices: G.~Gr\"atzer, E.\,T. Schmidt, and R.\,W. Quackenbush, Acta Sci. Math. \\(Szeged) 
\tbf{70} (2004), 473--494.

\end{slide} 

\end{document}